\documentclass[11pt]{article}
\usepackage[all]{xy}
\usepackage{graphicx}
\usepackage[margin=2.95cm]{geometry}
\usepackage{amssymb}
\usepackage{amsmath, amsthm}
\usepackage{latexsym} 
\newtheorem{thm}{Theorem}[section]

\newtheorem{prop}[thm]{Proposition}

 \newtheorem{rmk}[thm]{Remark}
\newtheorem{ex}[thm]{Example} 

\newcommand{\pf}{\noindent{\bf Proof.}\ }
\newcommand{\complex}{{\mathbb C}}

\newcommand{\reals}{{\mathbb R}}

\newcommand{\integers}{{\mathbb Z}}

\newcommand{\cinf}{{C^\infty}}
\newcommand{\sgn}{{\rm sgn~}}

\newcommand{\Hom}{{\rm Hom}}

\renewcommand{\ker}{{\rm Ker\, }}

\newcommand{\vbar}{{\overline{V}}}

\newcommand{\calc}{{\cal C}}
\newcommand{\cald}{{\cal D}}
\newcommand{\cale}{{\cal E}}
\newcommand{\calf}{{\cal F}}

\newcommand{\calh}{{\cal H}}

\newcommand{\calm}{{\cal M}}
\newcommand{\caln}{{\cal N}}
\newcommand{\calo}{{\cal O}}

\newcommand{\cals}{{\cal S}}

\newcommand{\calu}{{\cal U}}
\newcommand{\calv}{{\cal V}}

\newcommand{\del}{\partial}

\newcommand{\half}{\textstyle{\frac{1}{2}}}

\newcommand{\arrows}{\,\lower1pt\hbox{$\longrightarrow$}\hskip-.24in\raise2pt
             \hbox{$\longrightarrow$}\,}

\begin{document}

\title{{\bf The Maslov cycle as a Legendre singularity and projection of a wavefront set}}
\author
{Alan
Weinstein
\thanks{Research partially supported by NSF Grant
DMS-0707137 and the France-Berkeley Fund.
\newline \mbox{~~~~}MSC2010 Subject Classification Number: 
53D12 (Primary), 58J40, 81S10 (Secondary).
\newline \mbox{~~~~}Keywords: symplectic vector space, lagrangian
grassmannian, Fourier integral distribution, Maslov cycle}
\\
Department of Mathematics\\ University of California\\Berkeley, CA
94720 USA\\ {\small(alanw@math.berkeley.edu)}}
\maketitle
\begin{abstract} 
A Maslov cycle is a singular variety in the lagrangian grassmannian
 $\Lambda(V)$ of a symplectic vector space $V$
consisting of all lagrangian subspaces
 having nonzero intersection with
a fixed one.   Givental has shown that a Maslov cycle
 is a Legendre singularity, i.e. the projection of a smooth conic lagrangian submanifold $\cals$ in the cotangent bundle of $\Lambda(V)$.  We show here that $\cals$ is the wavefront set of a Fourier integral distribution which is
``evaluation at $0$ of the quantizations."
\end{abstract}

\section{Introduction}
\label{sec-intro}

Let $V$ be a symplectic vector space of dimension $2n$ over $\reals$ , 
$\Lambda=\Lambda(V)$ its
grassmannian of lagrangian subspaces.   
We refer to \cite{du:fourier} for more details concerning the following
basic facts concerning the lagrangian grassmannian.

Fix some $L_0 \in \Lambda$ and define the {\bf Maslov
  cycle} $\Sigma$ (with respect to $L_0$) 
to be the subset of $\Lambda$ consisting of those elements
which have a nonzero intersection with $L_0$.   $\Sigma$ is the
union of smooth strata $\Sigma_k = \{L\in \Lambda | \dim (L\cap L_0)
= k\}$, $k = 1,\ldots,n$.  
 $\Sigma_k$ is a locally
closed submanifold of $\Lambda$;  its closure is
the union of the $\Sigma_j$ for $j\geq k$.  In particular,
$\Sigma$ is the closure of $\Sigma_1$.   There is a natural coorientation of $\Sigma_1$; since the singularities of $\Sigma$ have codimension at least $3$, $\Sigma$ is dual to a well-defined cohomology class $\mu_V \in H^1(\Lambda,\integers)$.  This class, which is independent of the choice of $\Lambda$, is called the {\bf Maslov class}.    Associated with this cohomology class is a complex line bundle $\calm_V$ over $\Lambda$ called the Maslov line bundle.   It is the flat line bundle whose holonomy is given by the action of $\integers$ on $\complex$ taking $k$ to multiplication by $i^k$.    The structure group of this flat bundle is thus $\integers_4$.

The tangent space to $\Lambda$ at each element $L$ may be identified
with the space $S^2(L)$ of symmetric bilinear forms on $L$, considered as maps from $L$ to $L^*$.   (There are several ``natural'' identifications, differing by signs or factors of 2.   In Section \ref{sec-tangent}, we will select one of them.)
Under this identification, 
the tangent space to $\Sigma_k$ at  $L$ consists of those
forms $A$ for which the kernel of the associated linear map
$A_\flat:L_0 \to L^*_0$ contains the $k$-dimensional subspace $L
\cap L_0$.   The normal space $N_L \Sigma_k = 
T_L \Lambda /T_L \Sigma_k$ may therefore
be identified with
$S^2(L\cap L_0)$; in particular,
the codimension of $\Sigma_k$ in $\Lambda$ is $k(k+1)/2$.  

When the cotangent
space $T^*_L \Lambda$ is identified with the symmetric tensor product
$S_2(L) = L \circledcirc L$ of $L$ with itself,
the conormal space  
$N^*_L \Sigma_k$ for $L\in \Sigma_k$
becomes the subspace $(L\cap L_0)\circledcirc (L\cap L_0)$.  
Like the punctured\footnote{If $E$ is any vector bundle, its {\bf punctured} version is the complement of the zero section; we denote this complement by $\dot{E}.$}  conormal bundle to any submanifold,
$\dot{N}^*\Sigma_k$ is a conic 
lagrangian submanifold of the punctured cotangent bundle $\dot{T}^*L$; like $\Sigma_k$ itself, it is
not closed.   

Inside the union of the punctured conormal bundles to the strata of
$\Sigma$
is the conic subset
$\cals = \{(L,-v \circledcirc v) \in \dot{T}^*\Lambda | v\in L\cap L_0,
v\neq 0 \}.$ 
In each fibre $T_L^*\Lambda$ over
a point of $\Sigma_k$, the linear span of the fibre of $\cals$ is
the entire conormal space to $\Sigma_k$.  

A result of Givental (\cite{gi:nonlinear}, Section 10), specialized from lagrangian submanifolds to lagrangian subspaces, shows that $\cals$ is a smooth, conic
lagrangian submanifold of 
$\dot{T}^*\Lambda$ which is the closure of its restriction to
$\Sigma_1$.   This restriction is the 
side of $N^* \Sigma_1$ which is commonly used to define the 
coorientation to the Maslov cycle.  (Compatibility with this coorientation is
our reason for choosing
$-v\circledcirc v$ rather than $v\circledcirc v$.)
We will prove that $\cals$ is smooth by exhibiting non-degenerate local phase
functions for it.

Given any smooth conic lagrangian submanifold in a cotangent
bundle $T^*M$,  its  projection to  $M$  is known as a {\bf Legendre singularity} or {\bf
  wave front}).  Thus, a Maslov cycle is a wave front.
Another example, given in \cite{gu-ka-sh:sheaf}, is that
of the light cone in $\reals^{n+1}$ for any $n$, which is the
projection of the characteristic relation of the wave equation.  
 It is not clear to us how 
special the Legendre singularities are among the
singular hypersurfaces.

The main result in this note is the realization of $\cals$ as the wavefront set of  a Fourier integral
distribution $\phi$ on the lagrangian grassmannian.  This
distribution is the extension of a locally $L^1$ function
given by the evaluation at $0$ of the quantizations of 
lagrangian subspaces.   More precisely, it is the extension of a locally $L^1$ section
of the dual of a line bundle $E$ consisting of such quantizations.

The work described in this paper is an incarnation of the speculative idea that certain ``impossible'' operations on distributions, such as multiplication and restriction to submanifolds, might be definable as distributions on spaces of distributions.  To define such a distribution, for instance a generalized function $\psi$ from $\cald$ to $\complex$ such a definition, the relevant structure on a space $\cald$ of distributions might be its diffeology, i.e. the collection of smooth maps $\theta : P\to \cald$, known as plaques, from finite dimensional parameter spaces $P$.  For suitable plaques, there should be a distribution $\psi_\theta$ on $P$ which plays the role of the pullback of $\psi$ by $\theta$ in the sense that, for smooth maps $\delta:P' \to P$ (perhaps just submersions),  the compatibility condition $\psi_{P'} = \delta^* (\psi_P)$ is satisfied.  We hope that the example worked out in this paper, the pullback of an ``evaluation distribution'' on the distributions on a vector space by a particular plaque,  will eventually become part of a broader theory.

{\bf Acknowledgments} I would like to thank the group in Analyse
Alg\'ebrique of the Institut
Math\'ematique de Jussieu for many years of hospitality, 
Pierre Schapira and Bernd Sturmfels for directing my attention to
their work related to the problem discussed here, Shamgar Gurevich for many hours of discussion about the quantization of lagrangian subspaces,  Alexander Givental for encouraging me to extend his work in the quantum direction, and Luke Oeding and Maciej Zworski for advice on determinants and Fourier transforms.

\section
{The cotangent bundle of the lagrangian grassmannian} 
\label{sec-tangent}

There are two common ways to identify the tangent space to $\Lambda=\Lambda(V)$ at a point $L$
 with the space $S^2(L)$ of symmetric bilinear forms
on $L$.
Unfortunately, they differ by a sign.

First, we may consider $\Lambda$ as a homogeneous space for the
linear symplectic group $\mathrm{Sp}(V)$, whose Lie algebra $\mathfrak{sp}(V)$ may be
identified with the space $S^2(V)$ of symmetric bilinear forms on $V$.  This is done by identifying 
symmetric bilinear forms with quadratic functions and then with their
 (linear) hamiltonian
vector fields.   The tangent space $T_L \Lambda$ may then be
identified with the quotient of $S^2(V)$ by the subspace corresponding
to the isotropy subalgebra of $L$.  Since the quadratic functions whose
hamiltonian vector fields preserve $L$ are those which are constant,
hence zero, on $L$, the quotient in question may be identified
with the quadratic functions, or symmetric bilinear forms, on $L$
itself.

In Section \ref{sec-classandsymbol}, we will  interpret this identification in terms of the momentum map
for the action of $\mathrm{Sp}(V)$ on the cotangent bundle of $\Lambda$.

For the second identification, we choose a lagrangian subspace $N$
 transversal to $L$ in $V$.  The symplectic structure on $V$ gives an identification of $N$ with $L^*$ in such a way that $V$ is identified with $T^*L = L\times L^*$ with the canonical symplectic structure. 
Let $\calu_N$ denote the open subset of
$\Lambda$ consisting of those subspaces which are transversal to $N$.
Each $M \in \calu_N$ is the image of a closed linear one-form on $L$
which is the differential of a quadratic function corresponding to 
symmetric bilinear form $\rho_N(M)$ on $L$.    Equivalently, $M$ is the graph of the symmetric operator from $L$ to $L^*$ corresponding to $\rho_N(M)$.

The mapping $\rho_N$ is
a diffeomorphism taking $L$ to $0$ between $\calu_N$ and the vector space $S^2(L)$.  The differential $T_L\rho_N$ is thus an isomorphism
from $T_L \calu_N = T_L L$ to $T_0 S^2(L) = S^2(L)$.     We will see in a moment that this isomorphism is independent of $N$.

The two identifications above differ by a sign.   Look, for instance, at the case where
$V$ {\em is} $T^*L$, with $L$ identified with the zero section and $N$
with the fibre over $0$.  Given a quadratic function $H$ on $L$, we
may extend it to $T^*L$ by requiring it to be constant on fibres.  It
is easy to see that the hamiltonian flow of the extended function, at
time $t$, maps the zero section to the differential of the function
$-tH$ and not that of $H$.  (The minus sign here is essentially the one
in the Hamilton equations $dp_i/dt = -\del H/\del q_i$.)  
Since the first identification does not depend on $N$, neither does the second.

\begin{rmk}
{\em If $N$ and $N'$ are two transversals to $L$, the ``transition
  map'' $\rho_{N'} \rho_N^{-1}$ is a diffeomorphism between open subsets
  of $S^2(L)$ which fixes $0$.  Since the differential at $0$ of $\rho_N$ is independent of $N$, 
  the differential of the transition map at $0$ is the identity.
  This may also be seen from the following explicit formula for the map as a
  rational map on the space of symmetric matrices:
$$\rho_{N'}(M) = (I + \rho_N(M)B_{NN'})^{-1} \rho_{N}(M),$$
where $B_{NN'}$ is the symmetric operator (independent of $M$) obtained by considering $N'$ as the graph of a map from $N\approx L^*$ to $L$.   The formula is obtained by algebraic manipulation of the 
simply derived formula
$$\rho_{N'}(M)^{-1} = \rho_N(M)^{-1} + B_{NN'}.   $$
}
\end{rmk}

\begin{rmk}
{\em
The entire discussion in this section extends to the ``nonlinear''
case, where $V$ is any symplectic manifold and $\Lambda$ is the
infinite dimensional manifold\footnote{For the present purposes, one may consider $\Lambda$ as a
  Fr\'echet manifold or as a diffeological manifold.}
 of lagrangian submanifolds in $V$.
In this case, $T_L\Lambda$ may be identified in two ways with the space
$Z^1(L)$ of closed 1-forms on $L$.  The first is by extending closed
1-forms on $L$ to an open neighborhood in $V$ and seeing how $L$ moves
under the corresponding locally hamiltonian flows.  The second is by
identifying a neighborhood of $L$ in $V$ with a neighborhood of the zero
section in $T^*L$ and, hence, a neighborhood of the point $L$ in
$\Lambda$ with a neighborhood of $0$ in $Z^1(L)$.
 Once again, the two identifications differ by a sign, and the second one is independent of the choice of cotangent structure around $L$.}
\end{rmk}

In what follows, we will always use the second type of identification,
i.e. the one which uses the identification of a neighborhood of $L$ in
$V$ with a neighborhood of the zero section in $T^*L$.  We note that
this identification is {\em opposite} to the one used to
define the coorientation of the Maslov cycle used in the
construction of the Maslov cohomology class.\footnote{ 
 In \cite{we:connections}
we mentioned both identifications without realizing that
they differ by a sign.  (Fortunately, the discrepancy does
not have any effect on the results of that paper.)
}   

\section
{Phase functions for $\cals$}

We return now to the Maslov cycle $\Sigma$
associated to $L_0\in \Lambda$, and to the conic subset $\cals \subset
\dot{T}^*\Lambda$ lying above it.
To show that $\cals$ is a smooth lagrangian submanifold,
we will represent it locally by phase functions, using local
parametrizations of $\Lambda$ by spaces of symmetric bilinear forms.

For each lagrangian complement $N$ to $L_0$, we have the
diffeomorphism $\rho_N$ to $S^2(L_0)$ from the open subset $\calu_N$
of $\Lambda$.
The intersection  $\Sigma \cap \calu_N$ 
is mapped 
by $\rho_M$ onto the affine variety in $S^2(L_0)$ consisting of the
degenerate bilinear forms, with $\Sigma_k \cap \calu_N$ going to the
locally closed submanifold  $S^2_k(L_0)$ consisting of the forms of
nullity $k$.  If a basis for $L_0$ is chosen, $S^2(L_0)$ becomes
identified with the space of $n\times n$ symmetric matrices, and the
closure of $S^2_k(L_0)$ is the {\bf determinantal variety} defined by
the vanishing of all the $(n-k+1) \times (n-k+1)$ minors.

We identify the dual of $S^2(L_0)$ with $S_2(L_0) = L_0\circledcirc L_0$ and the  cotangent bundle    $T^*S^2(L_0)$
with $S^2(L_0)\times S_2(L_0)$.  
Under the cotangent lift of $\rho_N$,
$\cals \cap T^* \calu_N$ maps to the set $\cals_N$ of pairs
$(A,-v\circledcirc v)$ for which $v\neq 0$ and $A_\flat v=0$,
where $A_\flat$
is the linear map $L_0 \to L^*_0$ associated to $A$. 

\begin{rmk}
{\em
According to \cite{ni-ra-st:algebraic}, for any quadratic 
form $A$ 
of nullity $k$, the
full conormal space $N^* \Sigma_k$ at $A$
consists of {\em all}
 $B \in S_2(L_0)$ for which $X_\flat B^\sharp =0$, where $B^\sharp$ is
 the linear map $L^*_0 \to L_0$ associated to $B$.  This conormal space
 is
spanned by the subset of $\cals$ lying over $A$, which consists of the elements of the conormal space of the form 
$-v\circledcirc v$, i.e. the nonpositive rank 1 elements.
 The union of the conormals to the
strata is the Zariski closure in $T^*\Lambda$ of the conormal space
$N^*\Sigma_1$ to the top stratum.  
}
\end{rmk}

We will now show that $\cals_N$ is a smooth lagrangian submanifold
of the punctured cotangent bundle by analyzing the following
generating function.

Let $S:S^2(L_0) \times \dot{L}_0\to \reals$ be the 
``tautological'' function 
$S(A,v) = -\half A(v,v)$.  The differential of $S$ with respect to the second
variable, at the point $(A,v)$, is -$A_\flat (v)$, The
fibre-critical set $C_S$ where this differential
vanishes is thus the punctured kernel $  \dot{\ker}A_\flat =  \{(A,v)| v\neq 0, A_\flat v = 0\}$, and the differential of
$S$ at $(A,v)$ with respect to $A$ is
$-\half v\circledcirc v$.  

We must check that $S$ is a nondegenerate phase function.  The
differential with respect to $(A,v)$ of the fibre-differential $-A_\flat (v)$
is the map $(\delta A,\delta v) \to -A_\flat (\delta v) - (\delta A)_\flat
v.$  This differential is surjective because, $v$ being non-zero, the map can
achieve any value in $L_0^*$ with a suitable choice of $\delta
A$, and $\delta v = 0$.  

We now know that $C_S$ is a manifold and that, by the theory of
generating functions (see, e.g. \cite{du:fourier}), 
the map $(A,v)\mapsto (A,-\half v\circledcirc v)$
from $C_S$ to $\dot{T}^* S^2(L_0)$ whose image
is $\cals_N$ is a lagrangian immersion.  This map is clearly two-to-one, so it
is a double covering of its image, which is therefore an embedded
lagrangian submanifold. 

\begin{ex}
{\em
Let $L_0 =\reals$, with coordinate $x$.  The space of symmetric bilinear
forms is 1-dimensional, and the generating function is
$S(a,x)=-ax^2.$  The critical set is the line $a=0$, and
the map into the cotangent bundle is given by $(0,x)\mapsto (0,-x^2)$.   
}
\end{ex}

Since, for any lagrangian subspace $L$, there is another such subspace
which is transversal to both $L$ and $L_0$, we have completed the
proof of the following result.

\begin{thm}
Let $V$ be a symplectic vector space, $L_0$ an element of the
lagrangian grassmannian $\Lambda$, and $\Sigma$ the Maslov cycle
consisting of elements of $\Lambda$ having nonzero intersection with
$L_0$.  Then the subset $\cals$ of the punctured cotangent bundle
$\dot{T}^*\Lambda$ consisting of all pairs $(L,-v\circledcirc v)$ for
which $v$ is a nonzero vector in $L\cap L_0$ is a smooth, closed
lagrangian submanifold.
\end{thm}

We will give another proof  in Section \ref{sec-classandsymbol}.

\section{Quantization of $\cals$ to a Fourier integral distribution}
The oscillatory integral construction of H\"ormander \cite{ho:fourier} (closely related to earlier work by Maslov \cite{ma:theory}) associates to the conic lagrangian submanifold 
$\cals$ of the punctured cotangent bundle $\dot{T}^*\Lambda$ a class of distributions 
on $\Lambda$,  whose wavefront sets are contained in $\cals$.
Each of these {\bf Fourier integral distributions} has a principal symbol which is a section of the Maslov line bundle of $\cals$, tensored with an auxiliary line bundle.
In this section, we will exhibit a particular distribution whose definition is quite natural and whose
wavefront set is all of $\cals$.  In the next section, we will look at its principal symbol.

\subsection{Densities and geometric quantization}

Before proceeding any further, we recall some basic definitions and facts concerning densities.   For any vector space $R$ over $\reals$, and any real\footnote{Actually, $\alpha$ may even be a complex number, but we will not use this possibility.}  number $\alpha$, the space $|R|^\alpha$ of $\alpha$-densities on $R$ is defined as the 1-dimensional complex vector space of maps $s$ from the nonzero elements of the top exterior power $\bigwedge^{\mathrm top} R$ to $\complex$ which are positively homogeneous of degree $\alpha$, i.e. for which $s(cv)= |c|^\alpha s(v)$ for all $c>0$.   
 If $R$ is a vector {\em bundle},  $|R|^{\alpha}$ will denote the line bundle whose fibres are the spaces of $\alpha$-densities on the fibres of $R$.  Finally, if $M$ is any manifold,  $|M|^{\alpha} $ will be used to denote $|TM|^{\alpha}$; when $M$ is also a vector space, we will make it clear when necessary which of the two possible meanings of  $|M|^{\alpha}$ applies.

 For any $\alpha$ and $\beta$, there is a natural product between $\alpha$-densities and $\beta$-densities which yields $(\alpha+\beta)$-densities.  A $0$-density is just a complex number, so $\alpha$- and $(-\alpha)$-densities on a vector space are naturally in duality.   On a manifold, it is the $\alpha$- and $(1-\alpha)$-densities which can be paired to give a $1$-density (sometimes called simply a density), which can be integrated to give a complex number if it is compactly supported.   Finally, any exact sequence $0\to P'\to P \to P''\to 0$ induces an isomorphism between $|P'|^\alpha \otimes |P''|^\alpha$ and $|P|^\alpha$.

Given any symplectic vector space $V$ with a distinguished lagrangian subspace $L_0$, the geometric quantization construction produces a vector space consisting of the sections of the (unique up to isomorphism) prequantization line bundle which are covariant constant along the leaves of the polarization consisting of the affine subspaces parallel to $L_0$.
These sections may be considered as sections of a line bundle over $V/L_0$, and the choice of a lagrangian complement $Q$ to $L_0$ allows us to identify such sections with sections of $|Q|^{1/2}$ (with $Q$ considered here as a manifold).   The smooth sections of  $|Q|^{1/2}$ are contained in the space  $\cald '(Q,|Q|^{1/2})$ of distributional, complex-valued half-densities on $Q$, defined as the topological dual space of compactly supported smooth half-densities.       

Given a  lagrangian complement $Q$ to $L_0$,  we will identify $V$ with the cotangent bundle $T^*Q$ in such a way that $Q$ is identified with the zero section and $L_0$ with the fibre over $0\in Q$.

\subsection{Pure symplectic spinors}

Our Fourier integral distribution on $\Lambda$ itself will actually be a section of the dual of the line bundle $E\to \Lambda$ of {\em pure symplectic spinors} \cite{ko:symplecticspinors},
obtained by geometric quantization.   We now describe this bundle.

Let $\calh V$ be the Heisenberg Lie algebra which extends the symplectic vector space $(V,\omega)$, i.e. $V\oplus i\reals$ with the Lie algebra bracket for which the $i\reals$ summand is central and for which the bracket $[u,v]$ of two elements of $V$ is the element $-i\omega(u,v)$ of $i\reals$.   
We will use the {\bf Schr\"odinger representation} $\sigma$ of $\calh V$ on $\cald '(Q,|Q|^{1/2})$.   
Here, for each element $u$ of the summand $Q$ in $V$, $\sigma(u)$ is $-i$ times the operator of Lie derivative by $u$, considered as a constant vector field on $Q$.    For an element $u$ of $L_0\cong Q^*$,   $\sigma(u)$ is multiplication by the corresponding linear function on $Q$.

A canonical basis $(e_1,\ldots,e_n,f_1,\ldots, f_n)$ adapted to the splitting $V=Q\oplus L_0$, induces canonical coordinates $(x_1,\ldots,x_n,\xi_1,\ldots,\xi_n)$ on $V$.   In these coordinates, $\sigma(e_j)$ is (Lie derivative by) $-i\del/\del x_j$ and $\sigma(f_j)$ is multiplication by $x_j$.

For each lagrangian subspace $L \in \Lambda$, the operators $\sigma(u)$ for $u\in L$ form a commuting family of first-order linear differential operators whose common solution space $E_L$ is 1-dimensional.    These solution spaces form, in a natural way, the fibres of a complex line bundle $E$ over $\Lambda$.    The elements of $E_L$ are known as the {\bf pure symplectic spinors} associated with $L$.

When $L=Q$, the operators given by the Schr\"odinger representation are the Lie derivatives by constant vector fields, and the solution space $E_Q$ consists of the translation-invariant half-densities, i.e. the constant multiples of $\sqrt{|dx|}$, where $dx$ is shorthand for $dx_1 \wedge \ldots \wedge dx_n$.  More generally, when $L$ is transversal to $L_0$ and is therefore the graph of a symmetric operator $B$ from $Q$ to $L_0 \sim Q^*$, $E_L$ consists of the constant multiples of $e^{i<x,Bx>/2} \sqrt{|dx|}.$    When $L$ belongs to the Maslov cycle, the elements of $E_L$ are no longer smooth.   In fact, 
for any $L$, the elements of $E_L$ consist of distributions supported on the subspace of $Q$
annihilated by the linear functionals making up $L\cap L_0$.  For instance, when $L=L^0$, $E_L$ consists of 
multiplies of the Dirac delta function times $\sqrt{|dx|}$. 

\begin{rmk}
\label{rmk-cbms}
\em
{It is essentially the content of Lecture 10 in
\cite{we:lectures} that the bundle $E$ of pure symplectic spinors is naturally isomorphic to the tensor product of the Maslov line bundle $\calm_V$ with the line bundle of densities on the tautological vector bundle over $\Lambda$.   By the latter, we mean the bundle whose fibre over $L$ is equal to $|L|$.
 (The statement in the cited lecture is not quite right; it is claimed that $E$ is the Maslov bundle.   But the calculations there justify the description as the tensor product with the densities.)
}
\end{rmk}

\begin{rmk}
\label{rmk-kernels}
\emph{Suppose that $V'$ and $V''$ are  symplectic vector spaces with distinguished lagrangian subspaces $L_0'$ and $L_0''$ respectively, and $V=V'\times \vbar''$, where $\vbar''$ is $V''$ with its symplectic structure multiplied by $-1$.
 Then the pure symplectic spinors in the quantization of $V$ relative to $L=L_0' \times L_0 ''$ may be thought of as kernels of  operators between the quantizations of $V'$ and $V''$.   The composition of these operators becomes a kind of quantization of the composition of linear canonical relations.   This ``category'' (in quotes because composition is only partially defined due to issues concerning the domains of the operators) was studied in detail in \cite{gu-st:problems} and then in \cite{tu-za:category}, where the pure symplectic spinors were called Fresnel distributions.
}
\end{rmk}

\subsection{The evaluation distribution $\phi$}
We first consider the evaluation of pure symplectic spinors at $0\in Q$ as a function $\phi$ on the part of $E$ lying over the complement of the Maslov cycle; i.e. $\phi(\theta) = \theta(0).$    This function takes its values in the line $|Q|^{1/2} \cong
|(V/L_0)|^{1/2}$, where $Q$ and $V/L_0$ are considered here as vector spaces rather than manifolds. Since $\phi$ is linear on the fibres of $E$, we will actually view it as a (partially defined, so far) section of the bundle $E^* \otimes 
| (V/L_0)|^{1/2}$ over $\Lambda$.

\begin{thm}
 The section $\phi$ (extended in any  way across the Maslov cycle $\Sigma$, which has measure zero) is locally $L^1$
and hence extends to a distributional section of  $E^* \otimes 
| (V/L_0)|^{1/2}$ on $\Lambda$ which is singular everywhere on $\Sigma$.
\end{thm}

\pf
Since the statement is local, we will analyze $\phi$ by using a covering of $\Lambda$ by coordinate systems with local trivializations of $E$.  One of these coordinate neighborhoods is $\Sigma_0=\calu_{L_0}$, where $\phi$ is smooth.  Next, we look at the open set $\calu _Q$ consisting of the subspaces transversal to $Q$.  Using an idea already found in \cite{ma:theory},  we think of $Q\oplus L_0$ as the cotangent bundle of $L_0$ rather than of $Q$ and use the Fourier transform to exchange the roles of $L_0$ and $Q$.  This works because the Fourier transformation exchanges multiplication and differentiation operators.   More concretely, we coordinatize $\calu_Q$ and trivialize the part of   $E$ lying over it by the map $$\Psi:(A,c)\mapsto \left(\rho_Q^{-1}(A),\calf^{-1}\left(c e^{-i<A\xi,\xi>/2}\sqrt{|d\xi |}\right)\right),$$ where $\calf^{-1}$ is the inverse Fourier transform, defined on half-densities by
$$ \calf^{-1}\left(f(\xi)\sqrt {|d\xi|}\right) = (2\pi)^{-n/2} \left(\int f(\xi)e^{i<x,\xi>}  d\xi\right)\sqrt{|dx|} .$$

  In particular,
$$\calf^{-1}\left(c e^{-i<A\xi,\xi>/2}\sqrt{|d\xi |}\right)= 
c (2\pi)^{-n/2}\int e^{i(-<A\xi,\xi>/2 + <\xi,x>)}\sqrt{|dx|} 
.$$

 When $A$ is invertible, we get an explicit formula for $\Psi(A,c)$ by substituting $\eta=\xi-A^{-1}x$ in the integral above to get
$$\calf^{-1}\left(c e^{-iA<\xi,\xi>/2}\sqrt{|d\xi |}\right)= c(2\pi)^{-n/2}e^{-i<A^{-1}x,x>/2}\int e^{-i<A\eta,\eta>/2}d\eta \sqrt{|dx|}.$$
The complex gaussian integral in the last formula is equal to 
$$(2\pi)^{n/2}|\det A|^{-1/2} e^{-(\pi i/4) \sgn A},$$
where the determinant of $A$ is taken with respect to the chosen coordinates.  (This choice cancels out because of the presence of $\sqrt{|dx|}$.)
The final result is: 
$$\Psi(A,c)= \left( \rho_Q^{-1}(A),c |\det A|^{-1/2}e^{i(<A^{-1}x,x>/2- (\pi/4)\sgn A)}\sqrt{|dx|}\right).$$

We see immediately that $\Psi(A,c)$ is indeed a pure symplectic spinor associated with $L$, the graph of the mapping $A^{-1}$ from $Q$ to $L_0$.   Furthermore,
setting $x=0$ shows that, up to multiplication by a locally constant function of unit norm on the set of invertible $A$, 
$\phi(\Psi(A,c))$ is equal to $ c|\det A|^{-1/2}\sqrt{|dx|}$.   Thus, proving local integrability on $\calu_Q$ is reduced to verifying that 
$|\det A|^{-1/2}$ is a locally integrable function of the symmetric matrix $A$.  This is proved in the proposition below.

Although not every $L\in \Lambda$ belongs to one of the two subsets treated above, we can find for any $L$ a subspace $Q'$ which is transversal to both $L_0$ and $L$ and reduce the local integrability of $\phi$ near $L$ to the previous case.
\qed

\begin{prop}
The function $|\det A|^{-1/2}$ on the space of $n \times n$ symmetric matrices is 
locally integrable. 
\end{prop}

\pf
For $n=1$, this is just the local integrability of the function $f(a)=|a|^{-1/2}$.   Assuming that the statement is true for $k \times k$ symmetric matrices for all $k < n$, we will prove it for $n$.

We first prove local integrability away from the origin.    If $A_0 \neq 0$, its kernel has dimension $k$ for some $k < n$.   This allows us, for $A$ in a neighborhood $\caln$ of $A_0$, to split the spectrum of $A$ (with multiplicities) into $k$ ``small'' eigenvalues and $n-k$ eigenvalues which remain bounded away from zero, with corresponding spectral projections $\Pi^0_A$ and $\Pi^1_A$ depending smoothly on $A$.  
Choosing  bases for the image of $\Pi^0_A$ depending smoothly on $A$, we obtain a submersion $s$ from  $\caln$ to the space of $k \times k$ symmetric matrices and a factorization 
$\det A = g(A) \det s(A),$ where $g(A)$ is bounded away from zero.    The integrability of $|\det A|^{-1/2}$ over $\caln$ now follows from the induction hypothesis.

For local integrability near $A_0=0$, we use the degree $n$ homogeneity of the determinant function.   To integrate $|\det A|^{-1/2}$ over the unit ball (with respect to some norm), we first integrate over the region between the spheres of radius $1/2$ and $1$, which is possible by the result for nonzero $A_0$.   Each time we contract the region by a factor of $1/2$, the integral is multiplied by the same factor of $2^{-n^2/2}$.  
(This multiplier is the product of the factor $2^{n/2}$ by which the integrand is multiplied and the factor $2^{-n(n+1)/2}$ by which volume is multiplied.)
The integral over the entire ball is thus the sum of a convergent geometric series, so it is finite.
\qed

To see that $\phi$ is a Fourier integral distribution, we will write it locally as an oscillatory integral, starting with our formula for the inverse Fourier transform:
$$\calf^{-1}\left(c e^{-i<A\xi,\xi>/2}\sqrt{|d\xi |}\right)= 
c (2\pi)^{-n/2}\int e^{i(-<A\xi,\xi>/2 + \xi(x))}d\xi\sqrt{|dx|} 
.$$
Evaluation at $0$ gives
$$\phi(\Psi(A,c)) = c (2\pi)^{-n/2}\int e^{-i<A\xi,\xi>/2}d\xi \sqrt{|dx|}
.$$

We change variables to modified polar coordinates to make the phase function homogeneous of degree 1.  First write $\xi = r\theta$, where $r$ is nonnegative and $\theta$ is a unit vector.   Then we have
$$\phi(\Psi(A,c)) = c(2\pi)^{-n/2}\int e^{-ir^2<A\theta,\theta>/2} r^{n-1}dr d\theta \sqrt{|dx|},$$  where $d\theta$ is the volume element on the unit sphere.   Making the further substitution $s=r^2$, we have $dr= \frac{1}{2}s^{-1/2}ds$, and we get:
\begin{equation}
\phi(\Psi(A,c)) = (c/2) (2\pi)^{-n/2}\int e^{-is<A\theta,\theta>/2} s^{n/2-1}ds d\theta \sqrt{|dx|},
\end{equation}
which we rewrite as
\begin{equation}
\label{eq-oscint}
 \phi(\Psi(A,c)) = (c/2) (2\pi)^{-n/2}\int e^{-is<A\theta,\theta>/2} s^{-n/2}s^{n-1} ds d\theta \sqrt{|dx|} 
\end{equation}
to be compatible with the viewpoint that $(s,\theta)$ are now the polar coordinates.
We may now prove: 

\begin{prop}
$\phi$  is a Fourier integral distribution on $\Lambda$ belonging to H\"ormander's class
$I^m(\Lambda,\cals,E^*)$, where $m=-\frac{1}{4}n(n+1)/2 $.
\end{prop}

\pf
Following the development of Fourier integral distributions in Chapter 25 of \cite{ho:analysisiv}, especially Proposition 25.1.5, we will write $\phi$ locally, modulo $C^\infty$ functions, as an oscillatory integral of the form 
$$\int e^{iS(A,\lambda)}a(s,\lambda) d\lambda, $$
where $S$ is a nondegenerate phase function, positively homogeneous of degree $1$ in $\lambda$, generating the lagrangian submanifold $\Lambda$, and $a$ is an amplitude of some order $k$.  The order $m$ of the FID is then related to  $k$ by the formula $m=k+\frac{1}{4}(2N-d),$ where $N$ is the number of integration variables and $d$ is the dimension of the manifold on which the distribution is defined.   In the present case, $N=n$ and $d=n(n+1)/2$.)

The only problem with the representation \eqref{eq-oscint} is that the amplitude $s^{-n/2}$ is not supported away from the origin.   But this is easily remedied.   Let $\chi (s)$ be a cutoff function, equal to $1$ for $s$ in $[0,1]$ and to $0$ for $s$ in $[2,\infty]$.  Then we can break the integral into two terms, multiplying the amplitude first by $\chi$ and then by $1-\chi$.  The first term gives a $\cinf$ function of $(A,c)$, since the integrand in 
\eqref{eq-oscint} is now integrable and of compact support, while the integrand in the second term is now supported away from the origin and is homogeneous of degree $-n/2$ for large $s$.  

It remains  to check that the Fourier integral distribution $\phi$ defined by the Fourier transform formula is in fact equal to the distribution $\phi'$ given by the locally $L^1$ function defined on the regular set.   So far, it is not excluded that  they could differ by a distribution supported on the Maslov cycle.  

Applying Theorem 4.4.7 of \cite{du:fourier}, which relates the order of a Fourier integral distribution to its containment in Sobolev spaces, we find that the ``local Sobolev order'' of $\phi$ is $0$, in the sense that $\phi$ belongs to the Sobolev space $(H^r)_{\mathrm{loc}}$ for $r < 0$ and not for $r>0$.  

The local  Sobolev order of $\phi '$ is found to be $0$ by the same argument used to prove its locally $L^1$ property, based on its growth of order ${-1/2}$ in the distance from  the Maslov cycle.   This means that the singular part 
 $\phi - \phi '$ must also have local Sobolev order at least 0.  This is impossible for a nonzero distribution supported on a codimension 1 variety, which can have local Sobolev order at most $-1/2$, which is the local Sobolev order of a delta function supported along a hypersurface. 
\qed

\section{The Maslov class of $\cals$ and the principal symbol of $\phi$}
\label{sec-classandsymbol}

As a lagrangian submanifold of a cotangent bundle, $\cals$ carries
its own Maslov class  $\mu_\cals$.    In this section, we will describe the
cohomology group  $H^1(\cals,\integers)$ in which the class lies and  compute $\mu_\cals$.  The Maslov line bundle $\calm_\cals$ with $\integers_4$ holonomy given by the Maslov class enters in the description of the 
 line bundle $\cale\to\cals$ in which the principal symbol of $\phi$ takes its values; we will describe this bundle in detail, and we compute the symbol, up to a constant factor.

\subsection{The Maslov class of $\cals$}

When $n = 1$, $\cals$ is 
a half-line, so its cohomology is trivial, and $\mu_\cals =0$.    

To describe $\cals$ for $n\geq 2$, we note that its elements are parametrized by pairs $(L,\pm v)$, where $L$ is a lagrangian subspace
and $v$ is a nonzero vector in $L\cap L_0$, determined up to a sign.   But it is hard to describe the manifold of all these pairs in terms of the projection $(L,\pm v)\mapsto L$, because the topology of the fibre over  $L$ (even its dimension) depends on the dimension of $L\cap L_0$.  It is much better to project to $\pm v$ instead.   From this point of view,\footnote{This is a symplectic version of the so-called Bott-Samelson desingularization of a Schubert variety, as described, for instance, on page 118 of \cite{ma:symmetric}.}
 we see that an element of $\cals$ consists first of an element $\pm v$ of $\dot{L_0}/\integers_2$ and then a lagrangian subspace $L$ containing $v$.   The space of such $L$, given a choice of $v$, is in natural bijective correspondence with the lagrangian grassmannian of $<v>^\perp/<v>$, where $<v>$ is the line through $\pm v$.  Thus, $\cals$ is the $\Lambda(\reals^{2n-2})$ bundle over $\dot{L_0}/\integers_2 \simeq PL_0 \times \reals^+$ ($\simeq$ denoting homotopy equivalence)   associated to the principal  $\mathrm{Sp}(2n-2)$ bundle for which the associated $\reals^{2n-2}$ bundle has fibre $<v>^\perp/<v>$ over $\pm v$.
 
The cohomology group $H^1(\cals,\integers)$ is naturally isomorphic to $\Hom(\pi_1(\cals),\integers)$, and   the fundamental group $\pi_1(\cals)$ is almost completely determined by the homotopy exact sequence
$$\cdots\rightarrow  \pi_2(\reals P^{n-1}) \rightarrow \pi_1(\Lambda(\reals^{2n-2})) \rightarrow \pi_1(\cals) \rightarrow
\pi_1(\reals P^{n-1})\rightarrow \mathbf{1},
$$
where $\mathbf{1}$ is the trivial group.   

We know, further, that the projection $\cals \to \dot{V}/\integers_2$ admits a cross-section given by $\pm v  \mapsto (L_0,\pm v)$, so the exact sequence above splits.  The existence of this cross section also implies that the leftmost map above is trivial, and we get the split short exact sequence
$$\mathbf{1} \rightarrow \pi_1(\Lambda(\reals^{2n-2})) \rightarrow \pi_1(\cals) \rightarrow
\pi_1(\reals P^{n-1})\rightarrow \mathbf{1}.
$$ 
 For $n\geq 2$, the action of $\pi_1(\reals P^{n-1})\cong \integers_2$ on $ \pi_1(\Lambda(\reals^{2n-2})) \cong \integers$ is trivial, because the structure group $\mathrm{Sp}(2n-2)$ of the bundle is connected.    For $n\geq 3$, we conclude that $\pi_1(\cals) \cong \integers \times \integers_2$, and hence $H^1(\cals,\integers) \cong \integers$, with the inclusion of each fibre inducing an isomorphism on integral $H^1$.   For $n=2$, $H^1(\cals,\integers) \cong \integers \times \integers$, with the inclusion of each fibre inducing the projection on the first factor.   

Each fibre, consisting of the lagrangian subspaces in $V$ which contain $\pm v$ for some $v$, is contained in the full grassmannian $\Lambda$, and it is not hard to see that this inclusion induces an isomorphism on integral $H^1$.   In particular, the Maslov class  $\mu_V$ pulls back under the cotangent bundle projection to a generator of integral $H^1$ on each fibre, which is the pullback of a well-defined class in $H^1(\cals,\integers)$.    This class is just the pullback $p^*\mu_V,$ where $p$ denotes here the restriction to $\cals$ of the cotangent bundle projection.  Now we have:

\begin{thm}
\label{thm-maslovpullback}
The Maslov class $\mu_\cals$ of the lagrangian submanifold $\cals \subset T^*\Lambda$ is the pullback $p_\cals^*(\mu_V)$
of the Maslov class in $H^1(\Lambda,\integers)$ under the restricted cotangent bundle projection.
\end{thm}

\pf
We will use the invariance of Maslov classes under coisotropic reduction, as established 
in  \cite{vi:intersections}.    Let $\calc$ be the submanifold
$ \{(L,-v \circledcirc v) \in \dot{T}^*\Lambda |v\in L
v\neq 0 \}$  of  $\dot{T}^*\Lambda$ containing $\cals$, defined in the same way but without the condition imposed
by $L_0$.   This submanifold is coisotropic at points of $\cals$ because it contains the lagrangian submanifold $\cals$; by invariance under the action of the symplectic group of $V$, it is coisotropic everywhere.     

The characteristic leaves in $\calc$ are just the fibres of the projection  $(L,-v \circledcirc v)\mapsto \pm v$ to $\dot{V}/\integers_2$.  One way to see this is 
 to use actions of $\mathrm{Sp}(V)$.
In the dual space $\mathfrak{sp}(V)^*$, identified with $V \circledcirc V$, the elements $v\circledcirc v$ with $v\neq 0$
form a coadjoint orbit $\calo '$ (a minimal nilpotent orbit).   The momentum map of the natural action on $V$ takes $v$ to $v\circledcirc v$ and so identifies $\calo'$ with $\dot{V}/\integers_2$.  In fact, it will be more convenient to use instead the other minimal coadjoint orbit, $\calo = -\calo '$, which is the image of the momentum map for the action on $\vbar$. Thus $\calo$ is naturally identified with 
$\dot\vbar/\integers_2$.

Note that there is a natural identification between the lagrangian grassmannians of $V$ and $\vbar$, but the coorientations of their Maslov cycles, and hence the Maslov classes $\mu_V$ and $\mu_\vbar$, are opposite to one another under this identification.

 The momentum map  of the lift to $T^*\Lambda$  of the action of $\mathrm{Sp}(V)$ on $\Lambda$ assigns to each cotangent vector $(L,B)$  the element $B$ of $L\circledcirc L$ considered as an element of $V\circledcirc V$.   Our manifold $\calc$ is just $J^{-1}(\calo)$; since $J$ is a Poisson map, we see again that $\calc$ is coisotropic.   It also follows from the fact that $J$ is a Poisson map that the (connected) fibres of $J:\calc \to \calo$ are the characteristic leaves.   (The restriction of $J$ to $\calc$ is a canonical relation from all of $T^*\Lambda$ to $\calo$.)   

We now look at several subbundles of the symplectic vector bundle  $T_\cals T^*\Lambda$.   The first is the coisotropic subbundle $T_\cals \calc$ and its symplectic orthogonal $(T_\cals \calc)^\perp$.   Their quotient is naturally isomorphic to the pullback $J^* (T\calo) \cong J^*( TV/\integers_2).$   In addition, we have the vertical bundle $\calv_\cals$ and the bundle $T\cals$.  The latter two are lagrangian, with the first once transversal to $T_\cals \calc$ and the second one contained in it.   In each case, the symplectic reduction to $J^* (T\calo) \cong J^*( TV/\integers_2)$ is a smooth lagrangian subbundle.   The fibres of the lagrangian grassmannian bundle $\Lambda(J^*( TV/\integers_2))$ are all canonically isomorphic to the lagrangian grassmannian $\Lambda$ of $V$, and hence the reduced subbundles are described by mappings $\gamma_{\calv_\cals}$ and $\gamma_{T\cals}$ from $\cals$ to $\Lambda$.   Once we have determined these mappings, we will use the invariance result from Section 2 of \cite{vi:intersections}
 to conclude that the Maslov class of $\cals\subset T\Lambda$ is equal to the difference 
of pullbacks
\begin{equation}
\label{eq-diffmaslov}
 \gamma_{T\cals}^*(-\mu_\vbar)   -  \gamma_{\calv_\cals}^*(\mu_\vbar) .     
\end{equation}

To compute the $\gamma$'s, we use the notation $(\delta_L,\delta (-v\circledcirc v))$ for the general element of $T\calc$, where 
the variations $\delta L$ and  $-\delta(v\circledcirc v)$ are constrained by the condition that $L$ must continue to contain $v$ as $v$ varies (or, equivalently, that $v$ must continue to remain in $L$ as $L$ varies).
For tangent vectors to $\cals$, $v$ must vary within $L_0$, and this is independent of the basepoint of the tangent vector, so the map $\gamma_{T\cals}$ is simply the constant map $\cals\to\Lambda$ with value $L_0$.   On the other hand, for vertical tangent vectors to $T^*\Lambda$ at $(L,v)$ which are also contained in $\calc$, the variation $\delta L$ is zero, and so $v$ varies within $L$ itself, and hence $\gamma_{\calv\cals}$ is the restriction to $\cals$ of the cotangent bundle projection $p$.

Using these computations and Equation (\ref{eq-diffmaslov}), we find that 
$\mu_\cals = 0-p^*(\mu_{\vbar})=p^*(\mu_V)$, and the proof is complete.
\qed

\subsection{The principal symbol of $\phi$}

We first review how the symbol of a Fourier integral distribution is defined (and calculated).

Suppose that the FID $u$ is a generalized section of a vector bundle $F$ over a manifold $M$, associated with a conic lagrangian submanifold $\cals \subset T^*M$.   Over  
 each small enough open subset $\calu$ of $M$, $u$ is defined by integrating, along the fibres of a submersion $p:P\to \calu$, a section of the bundle $|(\ker Tp)^*| \otimes p^* F $ over $P$.  This section has the form $e^{iS} a $, where the nondegenerate phase function $S\in \cinf( P)$ is a generating function for $\cals \cap T^*\calu$, and the amplitude $a$ is a section of $|(\ker Tp)^*| \otimes p^* F$ with nice asymptotic properties relative to the conic structure of $\cals$.  

The fibre-differential $d_{\rm fib}S$  is a section of $(\ker T)^*$.    Along the zero set $Z \subset P$  of $d_{\rm fib}S$, the intrinsic derivative of $d_{\rm fib}S$ is a bundle map from  $TP$ to $(\ker Tp)^* $.   Nondegeneracy of $S$ means that $Z$ is a manifold and  this bundle map is surjective, with kernel $TZ$.
The resulting exact sequence induces an isomorphism along $Z$ between  
$|Z|
\otimes |(\ker Tp)^* |$ and $|P|$.   Also, $Tp$ gives a surjection from $TP$ to 
 $p^* TM$, leading to an isomorphism along $Z$ between $|\ker Tp| \otimes p^* |M|$ and $|P|$.   Combining these two isomorphisms  and using the natural duality between
$|\ker T p|$ and $ |(\ker T p)^* |$, we get an isomorphism between $ |(\ker Tp)^* |^2$ and $|Z|\otimes p^*|M|^{-1}$, and hence between 
$ |(\ker Tp)^* |$ and $|Z|^{1/2}\otimes p^*|M|^{-1/2}.$  These isomorphisms convert the amplitude $a$ of a FID with phase function $S$ into a section $\tilde{a}$  of $|Z|^{1/2}\otimes p^*|M|^{-1/2}\otimes F$ over $Z$.   (Note here the natural appearance of half-densities!)

To obtain the symbol of the FID, we transfer the section $\tilde{a}$ to the lagrangian submanifold $\cals$ by the embedding of $Z$ onto $\cals\subset T^*M$,  we replace $|Z|$ by $|\cals|$, and we multiply by a ``Maslov factor'' in order to make the symbol invariant under a change of phase function (which can include the insertion of a quadratic function in extra variables).   The final result is that  the principal symbol takes values in the bundle
$\cale = p^*F \otimes p^* |\Lambda|^{-1/2} \otimes |\cals|^{1/2} \otimes \calm_\cals$ over $\cals$, where $\calm_\cals$ is the Maslov bundle, i.e. the flat complex line bundle whose holonomy is the reduction of the Maslov class $\mu_\cals$ from $\integers$ to $\integers_4$.

For the evaluation distribution $\phi$, $M$ is the lagrangian grassmannian $\Lambda$, and
$F$ is the line bundle $E^* \otimes 
| (V/L_0)|^{1/2}$; furthermore, by Theorem \ref{thm-maslovpullback}, we can replace $\calm_\cals$ by $p^*\calm_V$.
Thus, we have
$$\cale = p^*(E^* \otimes | (V/L_0)|^{1/2}) \otimes p^* |\Lambda|^{-1/2} \otimes |\cals|^{1/2} \otimes p^*\calm_V.$$

Next, we will analyze $|\cals|^{1/2}$ by using the exact sequence obtained by differentiating the maps in the fibration $\cals$ over $\dot{L_0}/\integers_2$ whose fibre over $\pm v$ is the lagrangian grassmannian of $<v>^\perp/<v>$.    The sequence attached to a point $(L,-v\circledcirc v)$ is
$$0\to S^2(L/<v>)\to T_{(L,-v\circledcirc v)}\cals \to L_0 \to 0.$$   Note that the arrow to $L_0$ is defined up to sign, but this ambiguity will disappear in the next step. 

From the exact sequence above, we get a natural isomorphism between $|T_{L,-v \circledcirc v)}\cals|^{1/2}$ and $|S^2(L/<v>)|^{1/2}\otimes |L_0|^{1/2}$.  This turns the fibre $\cale_{(L,-<v\circledcirc v)}$ into 
$$E_L^* \otimes |(V/L_0)|^{1/2} \otimes |S^2(L)|^{-1/2} \otimes |S^2(L/<v>)|^{1/2}\otimes |L_0|^{1/2}\otimes (\calm_V)_L.$$   

Now the cokernel of the (injective) pullback map from $S^2(L/<v>)$ to $S^2(L)$ is naturally isomorphic to $L^*\otimes <v>^*$.  (The information lost in restricting a symmetric bilinear form $A$ on $L$ to a complement of $<v>$ and then pulling back this restriction to $L$ is contained in the pairing between $L$ and $<v>$ induced by $A$.)   This gives a natural isomorphism between $|S^2(L)|^{-1/2} \otimes |S^2(L/<v>)|^{1/2}$ and $|L\otimes <v>|^{1/2}$.   Using the identity $\det (P\otimes Q) = (\det P)^m (\det Q)^n$ for endomorphisms $P$ and $Q$ of spaces of dimension $n$ and $m$ respectively, we find that
$|L\otimes <v>|^{1/2}$ is naturally isomorphic to $|L|^{1/2}\otimes |<v>|^{n/2}.$

Using the isomorphisms of the previous paragraph, along with the natural trivialization of 
$ |(V/L_0)|^{1/2}\otimes |L_0|^{1/2}$ coming from the symplectic pairing of $V/L_0$ with $L_0$, we arrive at  a natural isomorphism $$\cale_{(L,-v\circledcirc v)} 
\cong E_L^* \otimes  |L|^{1/2}\otimes |<v>|^{n/2} \otimes (\calm_V)_L.$$

Finally, we use the description of $E$ from Remark \ref{rmk-cbms}, as the tensor product  of the Maslov bundle and the densities on the tautological bundle, to replace $E_L^*$ by $(\calm_V)_L^* \otimes |L|^{-1/2}$.    This gives the simple final result:
$$\cale_{(L,-v\circledcirc v)} \cong |<v>|^{n/2}.$$

We may take a closer look at the bundle over $\cals$ whose fibre over $(L,-v\circledcirc v)$ 
is $|<v>|^{n/2}$.  An element of this fibre is a complex valued function $\rho$ on $<v>$ which satifies the homogeneity property $\rho(aw) = |a|^{n/2}\rho(w)$.
  A section of $\cale$ is, therefore, a function of 
triples $(L,-v\circledcirc v,w)$, where $v$ and $w$ are collinear nonzero vectors in $L\cap L_0$, which satisfies the homogeneity property above in $w$.
The function $|w/v|^{n/2}$ is a natural cross section.    (Note that the ambiguity in the sign of $v$ is not a problem.)

We conclude that the bundle $\cale$ has a natural trivialization, so the symbol of $\phi$, a section of $\cale$, may be identified with a complex function on $\cals$.  This function is invariant under the symplectomorphisms which leave $L_0$ invariant.   As a consequence, this function of pairs $(L,-v\circledcirc v)$ is independent of $v$ and, as a function of $L$, is constant on each stratum of the Maslov cycle.   Because the stratum $\Sigma_1$ is sense, the function is constant everywhere.    Aside from the value of this constant, this constitutes a complete description of the symbol of $\phi$.

\end{document}